\documentclass[reqno,11pt]{amsart}
\usepackage{amsmath,amssymb,amsfonts,
}

\usepackage{diagrams}

\diagramstyle[centredisplay,dpi=600,
nohug]
\newif\ifhavegoodboldmath\havegoodboldmathtrue
\newarrow{Tto}----{->}
\newarrow{Tinc}{\ifhavegoodboldmath bold\fi hook}---{->}
\newarrow{Tincdash}{\ifhavegoodboldmath bold\fi hook}{dash}{}{dash}{->}
\newarrow{Tonto}----{->>}
\newarrow{Tdash}{}{dash}{}{dash}{->}

\theoremstyle{plain}
\newtheorem{thm}{Theorem}[section]
\newtheorem{lem}[thm]{Lemma}
\newtheorem{cor}[thm]{Corollary}
\newtheorem{propo}[thm]{Proposition}

\theoremstyle{definition}
\newtheorem{defn}[thm]{Definition}
\newtheorem{ex}[thm]{Example}
\newtheorem{parraf}[thm]{}
\newtheorem*{ack}{Acknowledgments}

\theoremstyle{remark}
\newtheorem*{rem}{Remark}

\numberwithin{equation}{thm}

\newcommand{\HA}{\widehat{A}}
\newcommand{\HB}{\widehat{B}}
\newcommand{\hf}{\widehat{f}}
\newcommand{\hd}{\widehat{d}}
\newcommand{\HM}{\widehat{M}}

\newcommand{\BA}{\mathbb A}
\newcommand{\BD}{\mathbb D}
\newcommand{\NN}{\mathbb N}

\newcommand{\ZZ}{\mathbb Z}

\newcommand{\fp}{\mathfrak p}
\newcommand{\fq}{\mathfrak q}

\newcommand{\FS}{\mathfrak S}
\newcommand{\FT}{\mathfrak T}
\newcommand{\FU}{\mathfrak U}
\newcommand{\FV}{\mathfrak V}

\newcommand{\FX}{\mathfrak X}
\newcommand{\FY}{\mathfrak Y}
\newcommand{\FZ}{\mathfrak Z}

\newcommand{\sch}{\mathsf {Sch}}

\newcommand{\sfn}{\mathsf {NFS}}
\newcommand{\scha}{\mathsf {Sch}_{\mathsf {af}}}

\newcommand{\sfna}{\sfn_{\mathsf {af}}}

\newcommand{\CC}{\mathcal C}

\newcommand{\CF}{\mathcal F}

\newcommand{\CI}{\mathcal I}
\newcommand{\CJ}{\mathcal J}
\newcommand{\CK}{\mathcal K}

\newcommand{\CO}{\mathcal O}

\newcommand{\dirlim}[1]{\begin{array}[t]{c} \mathrm{lim}\\[-7.5 pt]
 {\longrightarrow} \\[-7.5 pt] {\scriptstyle {#1}} \end{array}}
\newcommand{\invlim}[1]{\begin{array}[t]{c} \mathrm{lim}\\[-7.5 pt]
 {\longleftarrow} \\[-7.5 pt] {\scriptstyle {#1}} \end{array}}

\newcommand{\lto}{\longrightarrow}
\newcommand{\xto}{\xrightarrow}

\newcommand{\epi}{\twoheadrightarrow}

\newcommand{\inc}{\hookrightarrow}

\newcommand{\dimp}{\Leftrightarrow}

\newcommand{\tr}{\triangle}
\newcommand{\tc}{\widehat{\otimes}}
\newcommand{\om}{\widehat{\Omega}}

\DeclareMathOperator{\spec}{Spec}

\DeclareMathOperator{\spf}{Spf}
\DeclareMathOperator{\supp}{Supp}

\DeclareMathOperator{\Img}{Im}

\DeclareMathOperator{\Ker}{Ker}
\DeclareMathOperator{\Hom}{Hom}
\DeclareMathOperator{\Der}{Der}
\DeclareMathOperator{\Homcont}{Homcont}
\DeclareMathOperator{\Dercont}{Dercont}
\DeclareMathOperator{\ga}{\Gamma}

\DeclareMathOperator{\fD}{\mathfrak{D}}

\DeclareMathOperator{\coh}{\mathsf{Coh}}
\DeclareMathOperator{\Com}{\mathsf{Comp}}
\DeclareMathOperator{\com}{\text{-}\mathsf{comp}}
\DeclareMathOperator{\modu}{\text{-}\mathsf{mod}}
\DeclareMathOperator{\Modu}{\mathsf{Mod}}
\DeclareMathOperator{\alg}{\text{-}\mathsf{Alg}}

\newcommand{\ie}{{\it i.e.} }
\newcommand{\lc}{{\it loc.cit.}}

\begin{document}

\title[Smoothness and Jacobi criterion on formal schemes]{Infinitesimal lifting and Jacobi criterion for smoothness on formal schemes}

\author[L. Alonso]{Leovigildo Alonso Tarr\'{\i}o}
\address{Departamento de \'Alxebra\\
Facultade de Matem\'a\-ticas\\
Universidade de Santiago de Compostela\\
E-15782  Santiago de Compostela, SPAIN}
\email{leoalonso@usc.es}
\urladdr{http://web.usc.es/\~{}lalonso/}

\author[A. Jerem\'{\i}as]{Ana Jerem\'{\i}as L\'opez}
\address{Departamento de \'Alxebra\\
Facultade de Matem\'a\-ticas\\
Universidade de Santiago de Compostela\\
E-15782  Santiago de Compostela, SPAIN}
\email{jeremias@usc.es}

\author[M. P\'erez]{Marta P\'erez Rodr\'{\i}guez}
\address{Departamento de Matem\'a\-ticas\\
Escola Superior de En\-xe\-\~ne\-r\'{\i}a Inform\'atica\\
Campus de Ourense, Univ. de Vigo\\
E-32004 Ou\-ren\-se, Spain}
\email{martapr@uvigo.es}

\thanks{This work was partially supported by Spain's MCyT and E.U.'s
FEDER research project MTM2005-05754.}


\hyphenation{pseu-do}

\begin{abstract} This a first step to develop a theory of smooth, \'etale and unramified morphisms between noetherian formal schemes. Our main tool is the complete module of differentials, that is a coherent sheaf whenever the map of formal schemes is of pseudo finite type. Among our results we show that these infinitesimal properties of a map of usual schemes carry over into the completion with respect to suitable closed subsets. We characterize unramifiedness by the vanishing of the module of differentials. Also we see that a smooth morphism of noetherian formal schemes is flat and its module of differentials is locally free. The paper closes with a version of Zariski's Jacobian criterion.

\end{abstract}

\maketitle


\section*{Introduction}

\setcounter{equation}{0}
One of the great achievements of Grothendieck's point of view in algebraic geometry was the relationship between the classical notion of simple point and the notion of infinitesimal lifting. He proved that a point that is ``geometrically simple'', \ie such that keeps being simple after an extension of base field, can be characterized by the existence of a lifting from any subscheme defined by a square zero ideal of an affine scheme to the full scheme. In recent times formal schemes are getting increasing importance due to the variety of applications in which they are involved, to name a few, as algebraic models of rigid spaces \cite{ray}, in the study of cohomology of singular spaces \cite{ha2} or, more recently, in the context of stable homotopy \cite{st}. One feels the need of a greater progress of the basic fundamentals of the theory of formal schemes, so far reduced more or less to the last chapter of \cite{EGA1} and parts of \cite{EGA31}.
This paper intends to be the first in a series in which infinitesimal conditions on locally noetherian formal schemes are explored together with their applications to cohomology. In a subsequent paper we will give the local structure of smooth and \'etale maps of formal schemes. In this first installment, we develop a theory of smooth morphisms for noetherian formal schemes. \emph{Chemin faisant}, we also treat the other properties related to infinitesimal lifting, namely \'etale and unramified morphisms.

These topics have already been treated in the literature, albeit very scarcely. Smoothness is studied by Yekutieli under a special hypothesis, specifically, condition (ii) in \cite[Definition 2.1]{y} corresponds to a smooth map in which the base is an ordinary noetherian scheme, so smooth formal embeddings are examples of smooth maps of formal schemes. There was also Nayak's 1998 thesis whose results were eventually incorporated to the treatise \cite{LNS}. They work in the slightly more general context of essentially pseudo finite type maps  (\textit{cf.} [\lc, \S 2.1]). Our work has been developed mostly in parallel to this. As there is some overlapping between this and [\lc], we will point it out in the appropriate place when it arises. In fact, both groups of authors have reached an agreement on terminology and their definition of module of differentials [\lc, beggining of \S 2.6] agrees with ours when both are defined.

Let us discuss the contents of this paper.
The paper begins with some preliminaries to ease the task of the reader. They are collected into the first paragraph. In the second, we establish the notions that we will study. Our definition is taken from the one in \cite[\S17.3]{EGA44} for topological algebras. Therefore we will define formal smoothness for a map of formal schemes as the existence of liftings from a map of \emph{ordinary} schemes $T \inc Z$ given by a square zero Ideal. This agrees with the definition of formal smoothness for topological algebras and looks very much like the only reasonable convention. We therefore consider the maps like $T\inc Z$ as test morphisms for the condition of being formal smooth, unramified or \'etale. We obtain that maps of formal schemes $\FT \inc \FZ$ given by a square zero Ideal also detect formal smoothness (Proposition \ref{levantform}). Next we add the condition  of being of pseudo finite type to define the notions of smoooth, unramified and \'etale morphism. Our task is to show that these notions behave in a pleasant way, as in the case of usual schemes. The section closes with the general properties of these notions.

The next section is devoted to the study of the right notion of cotangent bundle for formal schemes. This is the sheaf of differentials that is obtained completing the usual module of differentials. It is our basic tool for studying more advanced properties of smoothness. The definition guarantees that the sheaf of differentials is coherent for a pseudo finite type map of formal schemes. Its basic characterizing property is that together with the canonical derivation it represents the functor that associates to a sheaf of \emph{complete} $\CO_{\FX}$-modules, the module of \emph{continuous} derivations. After explaining the functoriality of our construction, we show the analogous of the two fundamental exact sequences in this context.

Once one is equipped with the tool of the module of differentials, one is able to show further properties like the fact that smoothness, unramifiedness and being \'etale are properties local on the base and also on the source (Proposition \ref{condinflocal}). We show that this properties pass from a map of usual schemes to a completion. We characterize unramifiedness by the vanishing of the module of differentials (Proposition \ref{modifcero}). Also we see that  a smooth morphism is flat and its module of differentials is locally free. Next we discuss the splitting of the fundamental exact sequence when one of the maps is smooth and the chapter closes with Zariski's Jacobian criterion in this context (Corollary \ref{critjacob}).

\begin{ack}
We have benefited form conversations on these topics and on terminology with Joe Lipman, Suresh Nayak and Pramath Sastry. We also thank Jos\'e Antonio \'Alvarez for his useful remarks and the Mathematics Department of Purdue University for hospitality and support. 

The diagrams were typeset with Paul Taylor's \texttt{diagrams.tex}.
\end{ack}

\section{Preliminaries}\label{sec0}
We denote by $\sfn$ the category of locally noetherian formal schemes, by $\sfna$ the subcategory of affine noetherian formal schemes and by $\sch$ the category of schemes. 

We will begin by recalling briefly some basic definitions and results about locally noetherian formal schemes. Of course, for a complete treatment  we refer the reader to \cite[\S 10]{EGA1}. We will give some detailed examples of formal schemes, which we will refer along this exposition, like the affine formal scheme and the formal disc. Next we deal with finiteness conditions for morphisms in $\sfn$, which generalize the analogous properties  in $\sch$. In  the class of  adic morphisms we recall the notions of finite type morphisms, already defined in \cite[\S 10.13]{EGA1}. In the wider class of  non adic morphisms we will study morphisms of pseudo finite type  (introduced in \cite[p. 7]{AJL1}\footnote{Morphisms of pseudo finite type have been also introduced independently by Yekutieli in \cite{y} under the name ``formally finite type morphisms''}). Last we will recall from \cite[Chapter \textbf{0}]{EGA44} some basic properties  of the completed module of  differentials $\om^{1}_{A/B}$ associated to a continuous morphism $A \to B$ of adic rings.

\begin{parraf}  \label{equiv}
\cite[(10.2.2) and (10.4.6)]{EGA1}
The functors
\begin{equation*} 
A \leadsto \spf(A) \qquad \mathrm{and} \qquad
\FX \leadsto \ga(\FX, \CO_{\FX})
\end{equation*}
define a duality between the category of adic noetherian rings and $\sfna$ that generalizes the well-known relation between the categories of rings and affine schemes.\end{parraf}

\begin{parraf} \label{not}
Every locally noetherian formal scheme is a direct limit of usual schemes and every morphism in $\sfn$ is a direct limit of morphisms of schemes. More precisely:

\begin{enumerate}
\item \cite[(10.6.3), (10.6.4)]{EGA1}
Given $\FX$ in $\sfn$ and $\CJ \subset \CO_{\FX}$ an Ideal of  definition, for all $n \in \NN$, $X_{n}$ will denote the scheme $(\FX, \CO_{\FX}/\CJ^{n+1})$. 
Then $\FX$ is the direct limit in $\sfn$ of the diagram of noetherian schemes
 $\{X_{n}, i_{mn}\colon X_{m} \inc X_{n}, m \leq n\}_{n \in \NN}$. We will recall this data saying that $\FX$ it is expressed as 
\[\FX = \dirlim  {n \in \NN} X_{n}\] with respect to the Ideal of definition $\CJ$
and leave implicit that the schemes $\{X_{n}\}_{n \in \NN}$ are defined by the powers of $\CJ$.
\item \cite[(10.6.7), (10.6.8) and (10.6.9)]{EGA1}
If $f:\FX \to \FY$ is a morphism in $\sfn$, given $\CJ \subset \CO_{\FX}$ and $\CK \subset \CO_{\FY}$ Ideals of  definition such that  $f^{*}(\CK)\CO_{\FX} \subset \CJ$, for each $n \in \NN$, $f_{n}\colon X_{n}:= (\FX, \CO_{\FX}/\CJ^{n+1}) \to Y_{n}:=(\FY, \CO_{\FY}/\CK^{n+1})$ will be the morphism of  schemes induced by $f$.
The morphism $f$ is the direct limit of the system $\{f_{n}\}_{n \in \NN}$ associated  to the Ideals of  definition $\CJ \subset \CO_{\FX}$ and $\CK \subset \CO_{\FY}$ and we will write 
\[f = \dirlim  {n \in \NN} f_{n}\]
\end{enumerate}
Henceforth we will use systematically the above  notations.
\end{parraf}
\begin{parraf}  \label{defadic}
\cite[\S 10.12.]{EGA1} A morphism $f\colon\FX \to \FY$ in $\sfn$ is \emph{adic} (or simply \emph{$\FX$ is a $\FY$-adic} formal scheme) if there exists an Ideal of  definition $\CK$ of  $\FY$ such that $f^{*}(\CK)\CO_{\FX}$ is an Ideal of  definition of  $\FX$. Note that if there exists an Ideal of definition $\CK$ of $\FY$ such that $f^*(\CK)\CO_{\FX}$ is an Ideal of definition of $\FX$, then all Ideals of definition of $\FY$ share this property.

 If $f$ is adic and $\CK \subset \CO_{\FY}$ is an Ideal of definition, then the  diagrams of schemes associated to the Ideals $\CK$ and $f^*(\CK)\CO_{\FX}$
\begin{diagram}[height=2em,w=2em,labelstyle=\scriptstyle]
X_{m}       &  \rTto^{f_{m}} & Y_{m}  &\\
\uTinc      &                & \uTinc &\\
X_{n}       &  \rTto^{f_{n}} & Y_{n}  & \qquad (m \ge n \ge 0)\\
\end{diagram}
are cartesian.
 
Composition of adic morphisms is an adic morphism and the adic property is stable under base-change in $\sfn$.

\end{parraf}

\begin{parraf} \label{defnenc} 
\cite[\S 10.14.]{EGA1}
Let $\FX$ be in $ \sfn$. Given $\CI \subset \CO_{\FX}$, a coherent Ideal, $\FX' := \supp(\CO_{\FX}/\CI)$ is a closed subset and  $(\FX', (\CO_{\FX}/\CI)|_{\FX'})$ is a locally noetherian formal scheme . We will say that $\FX'$ is the   \emph{closed (formal) subscheme} of  $\FX$ defined by $\CI$.

\cite[(10.4.4)]{EGA1}
Given  $\FU \subset \FX$ open, it holds that $(\FU,\CO_{\FX}|_{\FU})$ is a noetherian formal scheme  and we say that $\FU$ is an \emph{open subscheme  of  $\FX$}. 

A morphism $f:\FZ \to \FX$  is a \emph{closed immersion} (\emph{open immersion}) if there exists $\FY\subset \FX$ closed (open, respectively) such that  $f$ factors as
\[
\FZ \xto{g} \FY \inc \FX
\]
where $g$ is an isomorphism. 

Closed and open inmersions are adic morphisms.
\end{parraf}


\begin{defn} \label{defmtf}
A morphism $f:\FX \to \FY$ in $\sfn$ is of \emph{pseudo finite type} if there exist $\CJ \subset \CO_{\FX}$ and  $\CK\subset \CO_{\FY}$ Ideals of  definition with  $f^{*}(\CK)\CO_{\FX} \subset \CJ$ and such that the  induced morphism of  schemes, $f_{0}: X_{0}  \to Y_{0}$ is of  finite type. If $f$ is of pseudo finite type and adic we say that $f$ is of \emph{finite type} in agreement with \cite[(10.13.1)]{EGA1}. 
\end{defn}

We have the following examples of morphisms in $\sfna$ provided by \ref{equiv}:

\begin{ex}
Let $A$ be a $J$-adic noetherian ring and  $\mathbf{T} = T_{1},T_{2},\ldots,T_{r}$ a finite number of indeterminates. 
\begin{enumerate}
\item 
The ring of restricted formal series  $A\{\mathbf{T}\}$ is a $J \cdot A\{\mathbf{T}\}$-adic noetherian ring (\emph{cf.} \cite[(\textbf{0}, 7.5.2)]{EGA1}). We call $\spf (A\{\mathbf{T}\})$ the \emph{affine formal $r$-space over $A$} or the \emph{affine formal space of  dimension $r$ over $A$} and we will denote it by
$\BA_{\spf (A)}^{r}$. It is a model of the \emph{closed} disk in rigid geometry, \emph{cfr.} \cite[\S 2.2]{he}. Note that $\spf (A\{\mathbf{T}\}) = \spf(A) \times \spec(\ZZ[\mathbf{T}])$ is the base change on formal schemes of the affine space $\spec(\ZZ[\mathbf{T}])$ over 
$\spec(\ZZ)$, that is why we adopt this terminology. The canonical projection
\[
\BA_{\spf (A)}^{r}   \to \spf(A)
\]
is of finite type.
\item  
The formal power series ring $A[[\mathbf{T}]]$ is a $(J \cdot A[[\mathbf{T}]]+ \langle\mathbf{T}\rangle \cdot A[[\mathbf{T}]])$-adic noetherian ring (\emph{cf.} \cite[Theorem 3.3 and Exercise 8.6]{ma2}). We define  the \emph{formal $r$-disc over $A$} or \emph{formal disc of dimension $r$ over $A$} as
$
\BD_{\spf (A)}^{r} = \spf (A[[\mathbf{T}]])
$. It is a model of the \emph{open} disk in rigid geometry, \emph{cfr.} \cite[\S 2.3]{he}. It has no counterpart on usual schemes, so a name relating it to rigid geometry is convenient.
The natural projection
\[
\BD_{\spf (A)}^{r}   \to \spf(A)
\]
is of pseudo finite type.

\item 
Given an ideal $I \subset A$,  the   closed immersion 
\[
\spf(A/I) \inc \spf(A)
\] 
is  a finite type morphism.
\item
Let $a \in A$ and denote by $A_{\{a\}}$ the completion of $A_{a}$ with respect to the ideal $J \cdot A_{a}$. The morphism $A \to A_{\{a\}}$ induces the canonical inclusion in $\sfna$
\[
\fD (a) \inc \spf(A).
\]
It is a  finite type morphism.

\item
Given  $X' =\spec(A/I)$ a closed  subscheme of $X=\spec(A)$, let $\HA$ be the completion of $A$ with respect to the $I$-adic topology. The \emph{morphism of  completion of  $X$ along  $X' $}, $\kappa: X_{/X' }=\spf(\HA) \to X$
is of pseudo finite type and is of finite type  only if $X$ and $X' $ have the same underlying topological space hence,  $X_{/X'} = X$ and $\kappa=1_{X}$.
\end{enumerate}
\end{ex}

\begin{propo} \label{ptf}
Let $f:\FX \to \FY$ be in $\sfn$.
The morphism $f$ is of pseudo finite type  if, and only if,  for each $x \in \FX$, there exist  affine open subsets $\FV \subset \FY$ and $\FU \subset \FX$ with $x \in \FU$ and $f(\FU) \subset \FV$ such that $f|_{\FU}$ factors as
\[
\FU \overset{j} \to  \BD_{\BA_{\FV}^{r}}^{s} \xto{p} \FV
\]
where $r,\, s \in \NN$, $j$ is a  closed immersion and $p$ is  the canonical projection.  

If $f$ is of finite type, then the above factorization may be written, taking $s=0$, $\FU \overset{j} \to  \BA_{\FV}^{r} \xto{p} \FV$.

\end{propo}

\begin{proof}
Since this is a local property  we may assume $\FX = \spf(A)$ and $\FY = \spf(B)$. Given $J \subset A$ and $K \subset B$ ideals of  definition such that $KA \subset J$ let $f_{0}: X_{0}=\spec(A/J) \to Y_{0}= \spec(B/K)$ be the morphism induced by $f$.
As $f$ is  pseudo finite type, there exists a presentation 
\[\frac{B}{K} \inc \frac{B}{K} [T_{1},T_{2}, \ldots ,T_{r}] \overset{\varphi_{0}}\epi \frac{A}{J}.\] 
This  morphism lifts to a ring homomorphism 
\[B \inc  B [T_{1},T_{2}, \ldots ,T_{r}] \to A\]
that extends to a continuous morphism 
\begin{equation} \label{factotipofinit}
 B \inc  B\{\mathbf{T}\}[[\mathbf{Z}]]:=B \{T_{1},T_{2}, \ldots ,T_{r}\}[[Z_{1},Z_{2}, \ldots ,Z_{s}]] \xto{\varphi} A
\end{equation}
 such that the images of $Z_{i}$ in $A$ together with $KA$ generate $J$. Let $B' := B\{\mathbf{T}\}[[\mathbf{Z}]]$. It is easily seen that the  morphism of graded modules associated to $\varphi$
\[
\bigoplus_{n \in \NN} \frac{(KB'+ \langle\mathbf{Z}\rangle)^{n}}{(KB' +\langle\mathbf{Z}\rangle)^{n+1}} \xto{\mathrm{gr}(\varphi)}  \bigoplus_{n \in \NN} \frac{J^{n}}{J^{n+1}}
\]  
is surjective, therefore, $\varphi$  is also surjective (\cite[III, \S2.8, Corollary 2]{b}).

If $f$ is of finite type, we may take $K \subset B$ and $J \subset A$ ideals of definition such that $KA = J$, so we can choose $s=0$. Then, the factorization (\ref{factotipofinit}) may be written \[B \to  B \{T_{1},T_{2}, \ldots ,T_{r}\} \epi  A\] and corresponds with the one given in \cite[(10.13.1)]{EGA1}. 
\end{proof}

The next result is a general version of \cite[(10.3.5)]{EGA1} and follows from the corresponding property in $\sch$, \cite[(6.3.4)]{EGA1}.

\begin{propo}  \label{mtf}
We have the following:
\begin{enumerate}
\item  \label{mtf1}
Given $f\colon\FX \to \FY$ and $g\colon\FY \to \FS$ in $\sfn$,
if $f$ and $g$ are (pseudo) finite type morphisms, then  $g \circ f$ is  a   (pseudo) finite type morphism.
\item  \label{mtf3}
If $f\colon\FX \to \FY$ is a  (pseudo) finite type morphism, given $h\colon\FY' \to \FY$ a morphism  in $\sfn$ we have that $\FX \times_{\FY} \FY'$ is in $\sfn$ and that $f'\colon\FX_{\FY'} \to \FY'$ is of  (pseudo) finite type.
\item 
Take $\FX \overset{f}\to \FY \to \FS$ and $\FY \overset{g}\to \FS \to \FS$ in $\sfn$, such that $\FY \times_{\FS} \FY'$ is in $\sfn$. If $f$ and $g$ are (pseudo) finite type morphisms, then 
$f \times_{\FS} g \colon \FX \times_{\FS} \FX' \to \FY \times_{\FS} \FY'$
is of (pseudo) finite type.
\end{enumerate}
\end{propo}

\begin{proof}
By (\ref{defadic}) it suffices to prove  the assertions for pseudo finite type morphisms. First, (1) and (2) are deduced from  the corresponding \emph{sorites} in $\sch$.
Statement (2) follows from the formal argument in \cite[(\textbf{0}, 1.3.0)]{EGA1}. From this it follows that $\FX \times_{\FS} \FX'$ belongs to $\sfn$ as a consequence of (1) and (2). Now the result is a consequence of  the analogous property in $\sch$.
\end{proof}

The usual module of differentials of a homomorphism $\phi:A \to B$ of topological algebras is not necessarily complete, but its completion has the good properties of the module of  differentials in the discrete case.  

\begin{parraf}\label{unonueve}
(\emph{cf.} \cite[\S \textbf{0}, 20.4, p. 219]{EGA41})
Given $B \to A$  a continuous homomorphism  of  preadic rings\footnote{According to \cite [(\textbf{0}, 7.1.9)]{EGA1} a ring $A$ is \emph{preadic} if there exists an ideal of definition $J$ of $A$ such that the the collection $\{J^{n}\}_{n \in \NN}$ forms a fundamental system of neighborhoods of $0$ in $A$. If $A$ is moreover separated and complete then $A$ is \emph{adic}.} and $K \subset B,\, J \subset A$  ideals of  definition such that $KA \subset J$, we denote by $\om^{1}_{A/B}$, the completion of the $A$-module $\Omega^{1}_{A/B}$ with respect  to  the $J$-adic topology 
\[
\om^{1}_{A/B} = 
\invlim {n \in \NN} \frac{\Omega^{1}_{A/B}}{J^{n+1}\Omega^{1}_{A/B}}.
\]

The continuous $B$-derivation  $d_{A/B}: A \to\Omega^{1}_{A/B}$ extends  naturally, by Leibnitz' rule, to a continuous $\HB$-derivation which, with an abuse of  terminology,  we will call \emph{canonical complete derivation  of  $\HA$ over $\HB$}, and denote by
\[\hd_{A/B} \colon \HA \to \om^{1}_{A/B}.\]
The canonical complete derivation of $\HA$ over $\HB$ makes the diagram
\begin{diagram}[height=2em,w=2em,labelstyle=\scriptstyle]
A                    & \rTto^{d_{A/B}}   & \Omega^{1}_{A/B}\\
\dTto^{\mathit{can}} &                   & \dTto^{\mathit{can}} \\
\HA                  & \rTto^{\hd_{A/B}} & \om^{1}_{A/B}\\
\end{diagram}
commutative.

For each $n \in \NN$ let $A_{n} = A/J^{n+1}$ and $B_{n}/K^{n+1}$. There is a canonical identification
\[
\om^{1}_{A/B} \cong \invlim {n \in \NN} \Omega^{1}_{A_n/B_n}
\]
with which 
\[
\hd_{A/B} \cong \invlim {n \in \NN} d_{A_n/B_n}.
\]
\end{parraf}

\begin{rem}
Given $B \to A$ a morphism of  preadic rings, let $K \subset B,\, J \subset A$  ideals of  definition such that $KA \subset J$. As a consequence of the previous discussion there results that 
\begin{equation} \label{modifcomplto}
( \om^{1}_{A/B},  \hd_{A/B}) \cong  (\om^{1}_{\HA/\HB}, \hd_{\HA/\HB})
\end{equation}
where $\HA$ and $\HB$ denote the completions of  $A$ and $B$, with respect to the  $J$ and $K$-preadic topologies, respectively, and $\om^{1}_{A/B}$ and $\om^{1}_{\HA/\HB}$  denote the completions of $\Omega^{1}_{A/B}$ and $\Omega^{1}_{\HA/\HB}$, with respect to the $J$-preadic and $J\HA$-adic topologies, respectively.
\end{rem}

\begin{parraf}
Let $A \com$ be the category of  complete $A$-modules for the $J$-adic topology. For all $M \in A \com$ the isomorphism \[\Homcont_{A}(\Omega^{1}_{A/B},M)  \cong \Dercont_{B}(A , M) \quad \text{(\emph{cf.} \cite[(\textbf{0}, 20.4.8.2)]{EGA41})} \] 
induces the  following  canonical isomorphism of  $B$-modules
 \begin{equation} \label{derivmodulalgeb}
 \begin{array}{ccc}
 \Homcont_{A}(\om^{1}_{A/B},M) & \cong &\Dercont_{B}(\HA , M) \\
 u	&\rightsquigarrow& u \circ \hd_{A/B}.\\
 \end{array}
\end{equation}
In other words, the pair $(\om^{1}_{A/B}, \hd_{A/B})$ represents  the functor
\[
M  \in A \com \leadsto \Dercont_{B}(\HA,M).
\]
In particular, if $M$ is an  $A/J$-module we have the  isomorphism
\[
\Hom_{A}(\om^{1}_{A/B},M) \cong \Der_{B}(\HA , M).
\]
\end{parraf}


\begin{parraf} \cite[(10.10.1)]{EGA1} \label{deftriangulito}
Let $\FX=\spf(A)$ with $A$ a  $J$-adic noetherian ring, $X=\spec(A)$ and $X'=\spec(A/J)$, so we have that $\FX= X_{/X'}$. Given $M$  an $A$-module, $M^{\tr}$ denotes the topological $\CO_{\FX}$-Module
\[
M^{\tr}:=(\widetilde{M})_{/X'} = 
\invlim {n \in \NN} \frac{\widetilde{M}}{\widetilde{J}^{n+1}\widetilde{M}}.
\]
Moreover, a morphism $u: M \to N$ in $A\modu$ corresponds to a morphism of  $\CO_{X}$-Modules $\tilde u:\widetilde M \to \widetilde N$ that induces a morphism of  $\CO_{\FX}$-Modules 
\[M^{\tr} \xto{u^{\tr}} N^{\tr} \quad = \quad \invlim {n \in \NN} (\frac{\widetilde{M}}{\widetilde{J}^{n+1}\widetilde{M}} \xto{\tilde{u}_{n}} \frac{\widetilde{N}}{\widetilde{J}^{n+1}\widetilde{N}}).\]
So there is an additive covariant  functor from  the category of  $A$-modules to the category of  $\CO_{\FX}$-Modules
\begin{equation} \label{triangulito}
\begin{array}{ccc}
A\modu& \overset{\tr} \lto & \Modu(\FX)\\
M					&  \rightsquigarrow    &  M^{\tr}.\\
\end{array}
\end{equation}
\end{parraf}

\begin{parraf}\label{propitrian}
If $M \in A\modu$ and $\HM$ denotes the complete  module of  $M$ for the $J$-adic topology, from the definition of the functor $\tr$ it is easy to deduce that: 
\begin{enumerate}
\item (\emph{Cf.} \cite[proof of (10.10.2.1)]{EGA1}) $\ga(\FX,M^{\tr}) =\HM$.
\item For all $a \in A$, $\ga(\fD(a),M^{\tr})=M_{\{a\}}$.
\item \label{equivtrian}
\cite[(10.10.2)]{EGA1} The  functor $(-)^\tr$ defines an equivalence of  categories between finite type $A$-modules   and the category $\coh(\FX)$ of  coherent $\CO_{\FX}$-Modules. 
\item \label{trianexac}
\cite[(10.10.2.1)]{EGA1} The functor $(-)^\tr$  is exact on the category of  $A$-modules of  finite type.
\end{enumerate}
\end{parraf}

\begin{parraf}\label{hazdif} A consequence of the previous results is the following. Let us consider a morphism $f \colon \spf(A) \to \spf(B)$ in $\sfna$. Let  $a \in A$ and $b \in B$ such that $f(\spf(A)_{\{a\}}) \subset \spf(B)_{\{b\}}$, then
\[
(\om^{1}_{A/B})^\tr (\spf(A)_{\{a\}}) = {(\Omega^{1}_{A/B})}_{\{a\}} 
= \om^{1}_{A_{a}/B} = \om^{1}_{A_a/B_b} 
= \om^{1}_{A_{\{a\}}/B_{\{b\}}}.
\]
Therefore, the sheaf $(\om^{1}_{A/B})^\tr$ agrees with the presheaf defined on principal open subsets by $\spf(A)_{\{a\}} \leadsto \om^{1}_{A_{\{a\}}/B}$.
\end{parraf}

\section{Definitions of the infinitesimal lifting properties}\label{sec1}

In this section we extend Grothendieck's classical definition of infinitesimal lifting properties from the category  of  schemes (\emph{cf.}   \cite[(17.1.1)]{EGA44}) to the category of locally noetherian formal schemes and we present some of their basic properties. We will refer to a morphism of formal schemes $f \colon \FX \to \FY$ simply as a $\FY$-formal scheme if there is no risk of ambiguity. If $X = \FX$ is an ordinary scheme, we will say that $X$ is a $\FY$-scheme.

\begin{defn}\label{defcinf}
Let $f:\FX \to \FY$ be a morphism in $\sfn$. We say that $f$ is \emph{formally smooth (formally unramified or formally \'etale)} if it satisfies the following lifting condition:

\emph{For all affine $\FY$-scheme $Z$ and for each closed subscheme $T\inc Z$ given by a square zero Ideal $\CI \subset \CO_{Z}$ the induced map
\begin{equation} \label{condlevant1}
\Hom_{\FY}(Z,\FX) \lto \Hom_{\FY}(T,\FX)
\end{equation}
is surjective (injective  or bijective, respectively).}

So, $f$ is formally \'etale if, and only if,  is formally smooth and formally unramified.
\end{defn}

\begin{parraf}  \label{condinfalg}
Let $f: \spf(A) \to \spf(B)$  be in $\sfna$. Applying (\ref{equiv}), we obtain that $f$ is formally smooth (formally unramified or formally \'etale) if, and only if,   the topological $B$-algebra $A$ 
is  formally smooth (formally unramified or formally \'etale, respectively) (\emph{cf.}  \cite[(\textbf{0}, 19.3.1) and (\textbf{0}, 19.10.2)]{EGA41}).

The reference for basic properties of the infinitesimal lifting conditions on preadic rings is \cite[\textbf{0}, \S\S{} 19.3 and 19.10]{EGA41}.
\end{parraf}

Next proposition shows that the lifting condition (\ref{condlevant1}) extends to a wider class of test maps.  

\begin{propo} \label{levantform}
Let $f:\FX \to \FY$ be in $\sfn$. If $f$ is formally smooth (formally unramified or formally \'etale), then for all affine noetherian $\FY$-formal scheme $\FZ$ and for all closed formal subschemes $\FT \inc \FZ$ given by a square zero Ideal $\CI \subset \CO_{\FZ}$, the induced map
\begin{equation} \label{condlevant3}
\Hom_{\FY}(\FZ,\FX) \lto \Hom_{\FY}(\FT,\FX)
\end{equation}
is surjective (injective or bijective, respectively). 
\end{propo}

\begin{proof}
Let $\FT=\spf(C/I) \overset{j}{\inc} \FZ=\spf(C)$ be a closed formal subscheme given by a square zero ideal $I \subset C$ . Let $L \subset C$ be an ideal of  definition, writing $T_{n}=\spec({C}/({I+L^{n+1}}))$ and $Z_{n}=\spec({C}/{L^{n+1}})$, the embedding  $j \colon \FT \inc \FZ$ is expressed as (see \ref{not}.(2))
\[\dirlim {n \in \NN} (T_{n} \overset{j_{n}}\inc Z_{n}),\] 
where the morphisms $j_{n}$ are  closed immersions of  affine schemes defined by a square zero Ideal. Given $u: \FT \to \FX$ a $\FY$-morphism,  we will denote by $u'_{n}$  the morphisms $T_{n} \inc \FT \xto{u} \FX$ that make the diagrams
\begin{diagram}[height=2em,w=2em,labelstyle=\scriptstyle]
T_{n}          & \rTinc^{j_{n}} & Z_{n}\\	
\dTto^{u'_{n}} &                & \dTto\\
\FX            & \rTto^f        & \FY\\
\end{diagram} 
commutative for all $n \in \NN$. 

Suppose that $f$ is formally smooth. Translating the argument given in  \cite[(\textbf{0}, 19.3.10)]{EGA41} for topological algebras to the context of formal schemes we get a $\FY$-morphism 
\[v :=\dirlim {n \in \NN} (v'_{n}: Z_{n} \to \FX)\] 
that satisfies $v|_{\FT}=u$. The morphisms $\{v'_{n}\}_{n \in \NN}$ are constructed by induction and satisfy that $v'_{n}|_{T_{n}}=u'_{n}$ and $v'_{n}|_{Z_{n-1}}=v'_{n-1}$, for each $n >0$ (\emph{cf.}  \cite[(\textbf{0}, 19.3.10.1) and (\textbf{0}, 19.3.10.2)]{EGA41}).

If $f$ is formally unramified, assume there exist $\FY$-morphisms 
$v: \FZ \to \FX$ and $w: \FZ \to \FX$ such that  $v|_{\FT} = w|_{\FT} = u$. With  the notations established at the  beginning of  the proof  consider  
\[v = \dirlim {n \in \NN} (v'_{n}: Z_{n} \to \FX) \,\,\text{ and } \,\, w= \dirlim {n \in \NN} (w'_{n}: Z_{n} \to \FX)\]  such that the   diagram
\begin{diagram}[height=2em,w=2em,p=0.3em,labelstyle=\scriptstyle]
T_{n} &   \rTinc^{j_{n}}	    &  Z_{n}				       &  	   	         &      &\\
   & \rdTto_{u'_{n}}&\dTto^{v'_{n}}\dTto_{w'_{n}}& \rdTto 		&     &\\ 
   &			    &\FX	  			       & 	 \rTto^{f} & \FY&\\  
\end{diagram}
commutes.
By hypothesis we have that $v'_{n} = w'_{n}$, for all $n \in \NN$, and we conclude that 
\[v=\dirlim {n \in \NN}v'_{n}= \dirlim {n \in \NN}  w'_{n}=w.\qedhere\]
\end{proof}

\begin{parraf} \label{lgfefnr}
In Definition \ref{defcinf} the test morphisms for the lifting condition are closed subschemes of \emph{affine} $\FY$-schemes given by square-zero ideals. An easy patching argument gives that the uniqueness of lifting conditions holds for closed subschemes of \emph{arbitrary} $\FY$-schemes given by square-zero ideals (\cite[(17.1.2.(iv))]{EGA44}). This applies to formally unramified and formally \'etale morphisms.
\end{parraf}

\begin{cor} \label{corlevfor}  
Let $f:\FX \to \FY$ be in $\sfn$. If the  morphism $f$ is formally unramified  (or formally \'etale), then for all noetherian $\FY$-formal schemes $\FZ$ and for each closed formal subscheme $\FT \inc \FZ$  given by a square zero Ideal $\CI \subset \CO_{\FZ}$, the induced map 
\begin{equation} \label{condlevant4}
\Hom_{\FY}(\FZ,\FX) \lto \Hom_{\FY}(\FT,\FX)
\end{equation}
is injective (or bijective, respectively).
\end{cor}

\begin{proof}
Given $\{\FV_{\alpha}\}$ a covering of affine open formal subschemes of  $\FZ$, denote by $\{\FU_{\alpha}\}$ the  covering of affine open formal subschemes of  $\FT$ given by $\FU_{\alpha} = \FV_{\alpha} \cap \FT$, for all $\alpha$.  By \cite[(10.14.4)]{EGA1} $\FU_{\alpha} \inc \FV_{\alpha}$ is a  closed immersion in $\sfn$ determined by a square zero Ideal. Therefore, the proof follows the same line as \cite[(17.1.2.(iv))]{EGA44}.
\end{proof}

The study of infinitesimal properties in $\sch$ using the module of  differentials leads one to look at the class of finite type morphisms. Under this assumption there are nice characterizations of the infinitesimal lifting conditions in terms of the module of  differentials. We will consider two conditions for morphisms in $\sfn$ that generalize the property of being of finite type for morphisms  in $\sch$: morphisms of pseudo finite type (\cite[1.2.2]{AJL1}) and its adic counterpart, morphisms of  finite type (\cite[(10.13.3)]{EGA1}).

\begin{defn} \label{smunet}
Let $f:\FX \to \FY$ be in $\sfn$. The  morphism $f$ is \emph{smooth (unramified or \'etale)}  if, and only if, it is of  pseudo finite type and formally smooth (formally unramified or formally \'etale, respectively). If moreover $f$ is adic, we say that $f$ is \emph{adic smooth (adic unramified or adic \'etale, respectively)}. So $f$ is  adic smooth (adic unramified or adic \'etale) if it is of  finite type and formally smooth (formally unramified or formally \'etale, respectively).

If $f: X \to Y$ is in $\sch$, both definitions agree with the one given in \cite[(17.3.1)]{EGA44} and we say  that $f$ is smooth (unramified or \'etale, respectively).
\end{defn}

Using \ref{condinfalg} we will be able to describe a few basic examples of  morphisms in $\sfna$ that satisfy some of the infinitesimal lifting conditions (Example \ref{exconforinf}). Before all else, let us recall some of the properties of  the infinitesimal lifting conditions for preadic rings.
\begin{rem} \label{confinfhafin}
Let $B \to A$ be a continuous  morphism of  preadic rings and take $J \subset A,\, K \subset B$ ideals of  definition with $KA \subset J$. Given $J' \subset A,\, K' \subset B$ ideals such that  $K'A \subset J'$, $J \subset J'$ and $K \subset K'$,  if $A$ is a  formally smooth (formally unramified or formally \'etale) $B$-algebra for the  $J$ and $K$-adic topologies then,  we have  that $A$ is a  formally smooth (formally unramified or formally \'etale, respectively) $B$-algebra for the $J'$ and $K'$-adic topologies, respectively.
\end{rem}

\begin{lem} \label{condinfalgcompl}
Let $B \to A$ be a continuous  morphism of  preadic rings, $J \subset A$ and $K \subset B$ ideals of  definition with $KA \subset J$ and let us denote by $\HA$ and $\HB$ the respective completions of $A$ and $B$. The following conditions are equivalent:
\begin{enumerate}
\item
$A$ is a  formally smooth (formally unramified or formally \'etale) $B$-algebra
\item
$\HA$ is a  formally smooth (formally unramified or formally \'etale, respectively) $B$-algebra
\item
$\HA$ is a  formally smooth (formally unramified or formally \'etale, respectively) $\HB$-algebra
\end{enumerate} 
\end{lem}

\begin{proof}
It suffices to note that
\[
\Homcont_{B\alg}(A,C) \cong \Homcont_{B\alg}(\HA,C) \cong \Homcont_{\HB\alg}(\HA,C) 
\] 
for all discrete rings $C$ and all continuous  homomorphisms $B \to C$.
\end{proof}

\begin{ex} \label{exconforinf}
Put $\FX= \spf(A)$ with $A$ a $J$-adic  noetherian ring and let $\mathbf{T} = T_{1},\, T_{2},\,  \ldots,\,  T_{r}$ be a   finite number of indeterminates.
\begin{enumerate}
\item \label{restrformliso}
If we take in $A$ the discrete topology, from  the universal property of the polynomial ring it follows that  $A[\mathbf{T}]$ is a formally smooth $A$-algebra. Applying the previous remark and  Lemma \ref{condinfalgcompl} we have that the  restricted formal series ring  $A\{\mathbf{T}\}$ is a formally smooth $A$-algebra, therefore the  canonical morphism $\BA_{\FX}^{r} \to \FX$ is  adic smooth.

\item  \label{formformliso}
Analogously to the preceding example, we obtain that $A[[\mathbf{T}]]$ is a  formally smooth $A$-algebra, from which  we deduce that projection $\BD_{\FX}^{r} \to \FX$ is smooth.

\item  \label{allavesformetale}
If  we take in $A$ the discrete topology it is known that, given $a \in A$, $A_{a}$ is a  formally \'etale $A$-algebra. So, there results that the canonical inclusion $\fD(a)\inc \FX$ is  adic \'etale.

\item  \label{enccerradoformnoram}
Trivially, every surjective  morphism of  rings is formally unramified. Therefore, given an ideal $I \subset A$,  the   closed immersion $\spf(A/I) \inc \FX$ is  adic unramified.

\item  \label{compformetale}
If  $X'= \spec(A/I)$ is a closed  subscheme of  $X$,  $\kappa: X_{/ X'} \to X$, the  morphism of completion of  $X$ along  $X'$, corresponds through (\ref{equiv}) with the continuous morphism of  rings $A \to \HA$, where $\HA$ is the  completion of $A$ for the $I$-adic topology  and therefore, $\kappa$ is \'etale.
\end{enumerate}
\end{ex}

\begin{propo} \label{infsorit}
In the category $\sfn$ of  locally noetherian formal schemes  the following properties hold:
\begin{enumerate}
\item  \label{infsorit1}
Composition of smooth (unramified or \'etale)  morphisms   is a smooth  (unramified or \'etale, respectively) morphism.
\item  \label{infsorit2}
Smooth, unramified and \'etale character  is stable under base-change in $\sfn$.
\item  \label{infsorit3}
Product of  smooth (unramified or \'etale) morphisms is a smooth   (unramified or \'etale, respectively) morphism.
\end{enumerate}
\end{propo}

\begin{proof}
Keeping in mind that composition of pseudo finite type maps is a pseudo finite type map  and pseudo finite type character of a map is preserved under base-change (Proposition \ref{mtf}) the proof is similar to \cite[(17.1.3) (ii), (iii) and (iv)]{EGA44}
\end{proof}

\begin{propo}  \label{infsoritcor}
The assertions of the last proposition hold if  we change the infinitesimal conditions  by the corresponding infinitesimal adic conditions.
 \end{propo}
 
\begin{proof}
By the definition of  the infinitesimal adic conditions, it suffices to  apply the  last result and the \emph{sorites} of finite type morphisms (Proposition \ref{mtf}).
\end{proof}

\begin{ex}
Let $\FX$ be in $\sfn$ and $r \in \NN$. From Proposition \ref{infsorit}.(\ref{infsorit2}),  Example  \ref{exconforinf}.(\ref{restrformliso}) and    Example \ref{exconforinf}.(\ref{formformliso})
we get that:
\begin{enumerate}
\item
The morphism of  projection $\BA_{\FX}^{r}:= \FX \times_{\spec(\ZZ)} \BA_{\spec(\ZZ)}^{r} \to \FX$ is an adic  smooth morphism.
\item
The canonical  morphism $\BD_{\FX}^{r}:=\FX \times_{\spec(\ZZ)} \BD_{\spec(\ZZ)}^{r} \to \FX$ is smooth.
\end{enumerate}
\end{ex}

\begin{propo}\label{infencajes}
The following holds in  the category of  locally noetherian formal schemes:
\begin{enumerate}
\item \label{infencajes1}
A  closed immersion is adic unramified.
\item \label{infencajes2} 
An open immersion is adic \'etale.
\end{enumerate}
\end{propo}

\begin{proof}
Closed and open immersions (see \ref{defnenc}) are monomorphisms and therefore, unramified. On the other hand, open immersions are smooth morphisms.
\end{proof}

\begin{propo} \label{propibasinf}
Let $f: \FX \to \FY ,\, g: \FY \to \FS$ be two morphisms of  pseudo finite type in $\sfn$. 
\begin{enumerate}
\item
If $g \circ f$ is unramified, then so is $f$. 
\item
Let us suppose that $g$ is unramified. If $g \circ f$ is smooth   (or \'etale) then, $f$ is smooth   (or \'etale, respectively).
\end{enumerate}
\end{propo}

\begin{proof}
Item (1) is immediate. The proof of (2) is analogous to \cite[(17.1.4)]{EGA44}.
\end{proof}

\begin{cor} \label{corbasinf}
Let $f:\FX \to \FY$ be a  pseudo finite type morphism and  $g: \FY \to \FS$  an \'etale morphism. The  morphism $g \circ f$ is smooth  (or \'etale) if, and only if,  $f$ is smooth  (or \'etale, respectively).
\end{cor}

\begin{proof}
It is enough to apply  Proposition \ref{infsorit}.(\ref{infsorit1}) and Proposition \ref{propibasinf}. 
\end{proof}

\section{Differentials of a pseudo finite type map of formal schemes} \label{sec2}

Given $f: X \to Y$ a finite type morphism of  schemes, it is well-known that $\Omega^{1}_{X/Y}$, the module of  $1$-differentials of  $X$ over $Y$, is an essential tool to study the smooth, unramified or \'etale character of  $f$. In this section, we introduce the  module of  $1$-differentials for  a morphism $f:\FX \to \FY$ in $\sfn$ and discuss its fundamental properties, which will be used  in the characterizations of  the infinitesimal conditions  in Section \ref{sec3}. We cannot use the general definition for ringed spaces, because it does not take into account the topology in the structure sheaves.

The observation in \ref{hazdif} shows that the following definition makes sense.

\begin{defn} \label{tresuno}
Given $f: \FX \to \FY$ in $\sfn$ we call \emph{module of  $1$-differen\-tials of  $f$} or \emph{module of $1$-differentials of  $\FX$ over $\FY$} and we will denote it  by $\om^{1}_{f}$ or $\om^{1}_{\FX/\FY}$, the sheaf of  topological $\CO_{\FX}$-Modules locally given by $(\om^{1}_{A/B})^{\tr}$ (see \ref{hazdif}), for all open sets $\FU=\spf(A) \subset \FX$ and $\FV=\spf(B) \subset \FY$ with $f(\FU) \subset \FV$. 
Note that $\om^{1}_{\FX/\FY}$ has structure of   $\CO_{\FX}$-Module.

Let $f: \FX \to \FY $ be in $\sfn$ and $\CJ \subset \CO_{\FX}$ and $\CK \subset \CO_{\FY}$ be Ideals of  definition  such that $f^{*}(\CK)\CO_{\FX} \subset \CJ$. These Ideals provide us with an inverse system of derivations
\[
d_{\FX_n/\FY_n} \colon \frac{\CO_{\FX}}{\CI^{n+1}} \to {\Omega^{1}_{\FX_n/\FY_n}}, \quad 
n \in \NN.
\]
Let  $\hd_{\FX/\FY}: \CO_{\FX} \to \om^{1}_{\FX/\FY}$ be the morphism
\[
\invlim {n \in \NN} d_{X_{n}/Y_{n}} =
\invlim {n \in \NN} (\frac{\CO_{\FX}}{\CI^{n+1}} \xto{d_{\FX_n/\FY_n}}
{\Omega^{1}_{\FX_n/\FY_n}}).
\]
It is locally defined for all couple of  affine open sets $\FU=\spf(A) \subset \FX$ and $\FV=\spf(B) \subset \FY$ such that $f(\FU) \subset \FV$ by 
\(
\hd_{\FX/\FY}(\spf(A)) = \hd_{A/B} \colon A \to \om^{1}_{A/B}.
\)
This construction is independent of the Ideals of definition chosen for $\FX$ and $\FY$, see \ref{unonueve}.

The  morphism $\hd_{\FX/\FY}$ is  a continuous $\FY$-derivation and it is called the \emph{canonical derivation of  $\FX$ over $\FY$}. 
We will refer to   $(\om^{1}_{\FX/\FY}, \hd_{\FX/\FY})$ as the  \emph{differential pair of  $\FX$ over $\FY$}.
\end{defn}

\begin{parraf} \label{defmofif}
If  $X=\spec(A) \to Y=\spec(B)$ is a morphism of  usual schemes, there results that $(\om^{1}_{X/Y},\hd_{X/Y}) = (\Omega^{1}_{X/Y}, d_{X/Y})$ is the  differential pair of the morphism of  affine schemes  (\emph{cf.}  \cite[(16.5.3)]{EGA44}).
\end{parraf}

\begin{rem}
Our definition of the differential pair, $(\om^{1}_{\FX/\FY}, \hd_{\FX/\FY})$, of a morphism $\FX \to \FY$ in $\sfn$, agrees with the  one given in \cite[2.6]{LNS} where it is directly defined as
\[
\invlim {n \in \NN} (\frac{\CO_{\FX}}{\CI^{n+1}} \xto{d_{\FX_n/\FY_n}}
{\Omega^{1}_{\FX_n/\FY_n}}).
\]
\end{rem}

\begin{propo} \label{finitmodif} (\emph{cf.} \cite[Proposition 2.6.1]{LNS})
Let $f: \FX \to \FY $ be a morphism in $\sfn$ of  pseudo finite type. Then \(\om^{1}_{\FX/\FY}\) is a coherent sheaf.
\end{propo}

\begin{proof}
We may suppose that $f: \FX=\spf(A) \to \FY=\spf(B)$ is in $\sfna$. Let $J \subset A$ and $K \subset B$ be ideals of  definition such that $KA \subset J$.
By hypothesis we have that $B_{0}=B/K \to A_{0}=A/J$ is a  finite type morphism and therefore,  $\Omega^{1}_{A_{0}/B_{0}}$ is a finite type $A_{0}$-module. From \cite[(\textbf{0}, 20.7.15)]{EGA41})  it follows that $\om^{1}_{A/B}$ is a finite type $A$-module. Therefore, since $\om^{1}_{\FX/\FY} = (\om^{1}_{A/B})^{\tr}$ the result is deduced from  \ref{propitrian}.\ref{equivtrian}.
\end{proof}

Given $X \to Y$ a morphism of  schemes in \cite[(16.5.3)]{EGA44} it is established that $(\Omega^{1}_{X/Y},d_{X/Y})$ is the  universal pair of the  representable  functor $\CF \in  \Modu(X) \rightsquigarrow \Der_{Y}(\CO_{X}, \CF)$. In  Theorem \ref{modrepr} this result  is generalized   for a morphism $\FX \to \FY$ in $\sfn$. 

\begin{parraf}
Given $\FX $ in $\sfn$ and $\CJ \subset \CO_{\FX}$ an Ideal of  definition of  $\FX$ we will denote by  $\Com (\FX)$ the full subcategory of  $\CO_{\FX}$-Modules  $\CF$ such that 
\[\CF = \invlim {n \in \NN} (\CF \otimes_{\CO_{\FX}} \CO_{X_{n}}).\] 
It is easily seen  that the definition does not depend on the election of the Ideal of  definition of  $\FX$. 

For example:
\begin{enumerate}
\item
Given $\FX= \spf(A)$ in $\sfna$ and $J \subset A$ an ideal of  definition for all $A$-modules $M$, it holds that 
\[M^{\tr}= \invlim {n \in \NN} \frac{\widetilde{M}}{\widetilde{J}^{n+1}\widetilde{M}} \in \Com(\FX).\]
\item
Let $\FX$ be  in $\sfn$. For all  $\CF \in \coh (\FX)$, we have that 
\[\CF =\invlim {n \in \NN} (\CF \otimes_{\CO_{\FX}} \CO_{X_{n}})\] by \cite[(10.11.3)]{EGA1} and therefore, $\coh (\FX)$ is a full subcategory of  $\Com (\FX)$. Consequently, if $f:\FX \to \FY$ is a pseudo finite type morphism in $\sfn$, then $\om^{1}_{\FX/\FY} \in \Com(\FX)$ by (\ref{finitmodif}).
\end{enumerate}
\end{parraf}

Now we are ready to show that given $\FX \to \FY$ a morphism in $\sfn$, $(\om^{1}_{\FX/\FY}, \hd_{\FX/\FY})$ is the  universal pair for the representable functor
\[\CF \in \Com (\FX)       \leadsto     \Dercont_{\FY}(\CO_{\FX}, \CF).\]

\begin{thm} \label{modrepr} 
Let $f: \FX \to \FY $ be a morphism in $\sfn$. Then  the canonical map
\[
\begin{array}{ccc}
\Homcont_{\CO_{\FX}}(\om^{1}_{\FX/\FY},\CF) & \overset{\varphi}\lto & \Dercont_{\FY}(\CO_{\FX}, \CF)\\
u				  &\rightsquigarrow    &u \circ \hd_{\FX/\FY}\\
\end{array}
\] is an isomorphism for every $\CF \in \Com (\FX)$.
\end{thm}

\begin{proof}
It is a globalization of \cite[(\textbf{0}, 20.7.14.4)]{EGA41}. We leave the details to the reader.
\end{proof}

\begin{lem} \label{imagdirecompl}
Let $f: \FX \to \FY $ be a morphism in $\sfn$. If $\CF \in \Com (\FX)$ then
\[
f_{*} \CF= \invlim {n \in \NN} (f_{*} \CF \otimes_{\CO_{\FY}} \CO_{Y_{n}})
\]
and consequently, $f_{*} \CF$ is in $\Com(\FY)$.
\end{lem}

\begin{proof}
Let $\CJ \subset \CO_{\FX}$ and $\CK \subset \CO_{\FY}$ be  Ideals of  definition such that $f^{*}(\CK)\CO_{\FX} \subset \CJ$. For all $n \in \NN$ we have the canonical morphisms  $f_{*} \CF \to f_{*} \CF \otimes_{\CO_{\FY}} \CO_{Y_{n}}$ that induce the  morphism of  $\CO_{\FY}$-Modules
\[
f_{*} \CF \lto  \invlim {n \in \NN} (f_{*} \CF \otimes_{\CO_{\FY}} \CO_{Y_{n}})
\]
To see whether it is an isomorphism is a local question, therefore we may assume that $f = \spf(\phi) \colon\FX=\spf(A) \to \FY= \spf(B)$ is in $\sfna$, $\CJ= J^{\tr}$ and $\CK= K^{\tr}$ with $J \subset A$ and $K \subset B$ ideals of  definition such that $KA \subset J$. Then $M=\ga(\FY, f_{*} \CF)$ is a complete $B$-module for  the $\phi^{-1}(J)$-adic topology and since $K \subset \phi^{-1}(J)$ we have that
\[M=\invlim {n \in \NN} \frac{M}{K^{n+1}M}\]
and the  result follows.
\end{proof}

\begin{propo} \label{cbasemodif}
Given a  commutative diagram  in $\sfn$ of  pseudo finite type  morphisms
\begin{diagram}[height=2em,w=2em,labelstyle=\scriptstyle]
\FX	  	 &  \rTto 		& \FY\\
\uTto^{g}    &    			& \uTto\\ 
\FX' 		&    	  \rTto^{h}	&     \FY' \\
\end{diagram}
there exists a morphism of  $\CO_{\FX'}$-Modules 
$g^{*}\om^{1}_{\FX/\FY} \lto \om^{1}_{\FX'/\FY'}$
locally determined by $\hd_{\FX/\FY} (a) \otimes 1 \rightsquigarrow \hd_{\FX'/\FY'} g(a)$.
Moreover, if the  diagram is cartesian, the  above morphism is an isomorphism.
\end{propo}

\begin{proof}
The morphism \[\CO_{\FX} \to g_{*}\CO_{\FX'} \xto{g_{*} \hd_{\FX'/\FY'}} g_{*} \om^{1}_{\FX'/\FY'}\] is a continuous $\FY$-derivation. Applying Proposition \ref{finitmodif} and  Lemma \ref{imagdirecompl} we have that $g_{*} \om^{1}_{\FX'/\FY'} \in \Com(\FX)$ and therefore by  Theorem \ref{modrepr} there exists an unique morphism of  $\CO_{\FX}$-Modules $\om^{1}_{\FX/\FY} \to g_{*} \om^{1}_{\FX'/\FY'}$ such that  the  following diagram is commutative
\begin{diagram}[height=2em,w=2em,labelstyle=\scriptstyle]
\CO_{\FX} 		&    	  \rTto^{\hd_{\FX/\FY}}	&  \om^{1}_{\FX/\FY} \\
\dTto   &    			& \dTto\\ 
g_{*} \CO_{\FX'} 	 &  \rTto^{g_{*} \hd_{\FX/\FY}}			&g_{*} \om^{1}_{\FX'/\FY'}\\ 
\end{diagram}
Equivalently, there exists a morphism of  $\CO_{\FX'}$-Modules $g^{*}\om^{1}_{\FX/\FY} \to \om^{1}_{\FX'/\FY'}$ locally determined by $\hd_{\FX/\FY} a \otimes 1 \rightsquigarrow \hd_{\FX'/\FY'} g(a)$. 

Let us suppose that the  square of formal schemes in the statement of this proposition is cartesian. We may assume $\FX=\spf(A),\, \FY= \spf(B),\, \FY'=\spf(B')$ and $\FX' = \spf(A')$ with $A' = A \tc_{B} B'$. The induced topology in $\om_{A/B} \otimes_{A} A'$ is the one given by the topology of $A'$. As a  consequence of the canonical  isomorphism  of  $A'$-modules $\Omega^{1}_{A'/B'}  \cong \Omega^{1}_{A/B} \otimes_{A} A' $ (\emph{cf.}  \cite[(\textbf{0}, 20.5.5)]{EGA41}) it holds  that 
\[
\om^{1}_{A'/B'} \cong 
\om^{1}_{A \tc_{B} B'/B'} \underset{(\ref{unonueve})}\cong 
\Omega^{1}_{A/B} \tc_{A} A' \underset{\textrm{\cite[(\textbf{0}, 7.7.1)]{EGA1}}}\cong 
\om^{1}_{A/B} \tc_{A} A'.
\] 
Finally $ \om^{1}_{A/B} $ is an $A$-module of  finite type (see Proposition \ref{finitmodif}) hence, $\om^{1}_{A'/B'} \cong \om^{1}_{A/B} \otimes_{A} A'$.  
\end{proof}
With  the previous notations, if $\FY=\FY'$ the  morphism $g^{*}\om^{1}_{\FX/\FY} \lto \om^{1}_{\FX'/\FY}$ is denoted by $dg$ and is called \emph{the differential of  $g$ over $\FY$}.

\begin{cor}  \label{modifesqu}
Given $f:\FX \to \FY$ a  \emph{finite type} morphism in $\sfn$ consider $\CK\subset \CO_{\FY}$ and $\CJ =f^{*}(\CK)\CO_{\FX} \subset \CO_{\FX}$ Ideals of  definition that let us express  
\[f = \dirlim {n \in \NN} (f_{n}: X_{n} \to Y_{n}).
\]
Then
\[
\Omega^{1}_{X_{n}/Y_{n}} \cong \om^{1}_{\FX/\FY} \otimes_{\CO_{\FX}} \CO_{X_{n}}.
\] for all $n \in \NN$.
\end{cor}
\begin{proof}
Since $f$ is an adic  morphism, the diagrams  
\begin{diagram}[height=2em,w=2em,labelstyle=\scriptstyle]
\FX	&	\rTto^{f} &	\FY \\
\uTto			&	&		\uTto\\
X_{n}		&	\rTto^{f_{n}} &	Y_{n}\\
\end{diagram}
are cartesian, $\forall  n \ge 0$. Then the  corollary follows.
\end{proof}

In  the  following example we show that if the  morphism is not adic, the last  corollary does not hold.
\begin{ex} \label{contrmod0}
Let  $K$ be a field  and $p:\BD^{1}_{K} \to \spec(K)$ the  projection  morphism of  the formal disc of  $1$ dimension  over $\spec(K)$. Given the  ideal of  definition $\langle T \rangle \subset K[[T]]$ such that 
\[p = \dirlim {n \in \NN} p_{n}\] 
we have that $\Omega^{1}_{p_{0}}=0$ but, 
\[
\om^{1}_{p} \otimes_{\CO_{\BD^{1}_{K}}} \CO_{\spec(K)} = 
(\om^{1}_{K[[T]]/K})^{\tr} \otimes_{K[[T]]^{\tr}} \widetilde{K} \cong 
\widetilde{K} \neq 0.
\]
\end{ex}

We extend the usual First and Second Fundamental Sequences to our construction of differentials of pseudo finite type morphisms between formal schemes. They will provide a basic tool for applying it to the study of the infinitesimal lifting. Also, we will give a local computation based on the Second Fundamental Exact Sequence.

\begin{propo}{(First Fundamental Exact Sequence)}\label{primersef}
Let $f\colon\FX \to \FY$  and $g\colon \FY \to \FS$ be two morphisms in $\sfn$ of  pseudo finite type. There exists an exact sequence  of  coherent $\CO_{\FX}$-Modules 
\begin{equation} \label{sef1}
f^{*}\om^{1}_{\FY/\FS} \xto{\Phi} \om^{1}_{\FX/\FS} \xto{\Psi} \om^{1}_{\FX/\FY} \to 0\\
\end{equation}
where $\Phi$ and $\Psi$  are locally defined  by
\[
\hd_{\FY/\FS} b \otimes 1 \leadsto \hd_{\FX/\FS} f(b) \qquad 
\hd_{\FX/\FS} a Ê         \leadsto \hd_{\FX/\FY} a\\
\]
\end{propo}

\begin{proof}
This is a globalization of \cite[Lemma 2.5.2]{LNS} (see also \cite[(\textbf{0}, 20.7.17.3)]{EGA41}).
The morphism $\Phi\colon f^{*}\om^{1}_{\FY/\FS}  \to \om^{1}_{\FX/\FS}$ of  $\CO_{\FX}$-Modules is $df$, the differential of  $f$ over $\FS$. Since $\hd_{\FX/\FY}\colon \CO_{\FX} \to \om^{1}_{\FX/\FY}$ is a continuous $\FS$-derivation, from  Theorem \ref{modrepr} there exists a unique  morphism of  $\CO_{\FX}$-Modules $\Psi\colon \om^{1}_{\FX/\FS} \to \om^{1}_{\FX/\FY}$ such that  $\Psi \circ \hd_{\FX/\FS} = \hd_{\FX/\FY}$.

As for proving  the exactness we can reduce to the affine case and then it is the first part  of \cite[Lemma 2.5.5]{LNS} (see also \cite[(\textbf{0}, 20.7.17)]{EGA41}).
\end{proof}

\begin{parraf}
Let $f\colon\FX \to \FY$ be  a pseudo finite type morphism in $\sfn$, and $\CJ \subset \CO_{\FX}$ and $\CK \subset \CO_{\FY}$ be Ideals of  definition with $f^{*}(\CK) \CO_{\FX} \subset \CJ$ and \[f\colon \FX \to \FY = \dirlim {n \in \NN} (f_{n}\colon X_{n} \to Y_{n})\]
the relevant expression for $f$. For all $n \in \NN$, from  the First Fundamental Exact Sequence  (\ref{sef1}) associated to $X_{n} \overset{f_{n}}\to Y_{n} \inc \FY$, we deduce that $\om^{1}_{X_{n}/\FY} = \Omega^{1}_{X_{n}/\FY}= \Omega^{1}_{X_{n}/Y_{n}}$.
\end{parraf}

\begin{parraf}
Given $\FX' \overset{i} \inc \FX$ a  closed immersion in $\sfn$  we have that the  morphism $i^{\sharp}: i^{-1}(\CO_{\FX})  \to \CO_{\FX'}$ is an epimorphism. If $\CK:= \ker(i^{\sharp})$ we call $\CC_{\FX'/\FX}:=\CK/ \CK^{2}$ \emph{the conormal sheaf of  $\FX'$ in $\FX$}. 
 
It is easily  shown that $\CC_{\FX'/\FX}$ satisfies the following properties:
\begin{enumerate}
\item
It is a coherent $\CO_{\FX'}$-module.
\item
If $\FX' \subset \FX$ is a closed subscheme given by a coherent  Ideal $\CI \subset \CO_{\FX}$, then $\CC_{\FX'/\FX}=i^{*}(\CI/\CI^{2})$.
\end{enumerate}
\end{parraf}

\begin{propo}{(Second Fundamental Exact Sequence)} \label{segunsef}
Let $f:\FX \to \FY$ be a pseudo finite type morphism in $\sfn$, and $\FX' \overset{i} \inc \FX$ a closed immersion. There exists an exact sequence of coherent $\CO_{\FX'}$-Modules
\begin{equation} \label{sef2}
\CC_{\FX'/\FX} \xto{\delta} i^{*}\om^{1}_{\FX/\FY} \xto{\Phi} \om^{1}_{\FX'/\FY} \to 0
\end{equation}
\end{propo}
\begin{proof}
Morphism $\Phi$ is the differential of $i$ and is defined by 
$\hd_{\FX/\FY} a \otimes 1 \rightsquigarrow \hd_{\FX'/\FY} i(a)$
(Proposition \ref{cbasemodif}). If $\CI \subset \CO_{\FX}$ is the  Ideal that defines the closed subscheme $i(\FX') \subset \FX$ the  morphism $\delta$ is the one induced by $\hd_{\FX/\FY}|_{\CI}: \CI \to \om^{1}_{\FX/\FY}$. Again the exactness is consequence of \cite[(\textbf{0}, 20.7.20)]{EGA41}).
\end{proof}

As it happens in $\sch$ the Second Fundamental Exact Sequence  leads to a local  description of the  module of  differentials of  a pseudo finite type morphism between locally noetherian formal schemes. 
\begin{parraf} \label{genermodifptf}
Let  $f: \FX=\spf(A) \to \FY=\spf(B)$ be a morphism in $\sfna$ of  pseudo finite type, then it  factors as  (see Proposition \ref{ptf})
\[
\FX= \spf(A) \overset{j} \inc \BD^{s}_{\BA^{r}_{\FY}}= \spf(B\{\mathbf{T}\}[[\mathbf{Z}]] ) \xto{p} \FY= \spf(B)
\]
where $r,\, s \in \NN$, $\mathbf{T} = T_{1},T_{2},\ldots,T_{r}$ and $\mathbf{Z} = Z_{1},Z_{2},\ldots,Z_{r}$ two sets of indeterminates, $p$ is the canonical projection and $j$ is a  closed immersion given by an Ideal $\CI= I^{\tr} \subset \CO_{\BD^{s}_{\BA^{r}_{\FY}}}$. Consider also a system of generators $I= \langle  P_{1},  \ldots,  P_{k}\rangle  \subset B\{\mathbf{T}\}[[\mathbf{Z}]]$.

The Second  Fundamental  Exact Sequence (\ref{sef2}) associated to $\FX \overset{j} \inc \BD^{s}_{\BA^{r}_{\FY}} \xto{p} \FY$ 
corresponds through the equivalence  between the category of finite type $A$-modules and $\coh(\FX)$ \ref{propitrian}.\ref{equivtrian}, to  the sequence 
\begin{equation} \label{segsucfunafdis}
\frac{I}{I^{2}} \xto{\delta} \om^{1}_{B\{\mathbf{T}\}[[\mathbf{Z}]]/B} \otimes_{B\{\mathbf{T}\}[[\mathbf{Z}]]}A\xto{\Phi} \om^{1}_{A/B} \to 0.
\end{equation}
Let us use the following abbreviation $\hd= \hd_{B\{\mathbf{T}\}[[\mathbf{Z}]] /B}$. Since 
\[\om^{1}_{B\{\mathbf{T}\}[[\mathbf{Z}]]/B}   \cong \om^{1}_{B[\mathbf{T},\mathbf{Z}]/B}    \cong  \bigoplus_{i=1}^{r}  B\{\mathbf{T}\}[[\mathbf{Z}]] \hd T_{i} \oplus \bigoplus_{j=1}^{s}  B\{\mathbf{T}\}[[\mathbf{Z}]] \hd Z_{j} 
\] then  
$\{ \hd T_{1},\, \hd  T_{2}, \ldots,\, \hd  T_{r},\, \hd Z_{1},\, \ldots,\, \hd Z_{s} \}$ is a basis of the free $B\{\mathbf{T}\}[[\mathbf{Z}]] $-module  $\om^{1}_{B\{\mathbf{T}\}[[\mathbf{Z}]] /B}$. Therefore, if $a_{1},\, a_{2}, \ldots,\, a_{r},\, a_{r+1},\, \ldots,\, a_{r+s} $ are the images of  $T_{1},\, T_{2},\, \ldots,\, T_{r},\, Z_{1},\, \ldots,\, Z_{s}$ in $A$, by the definition of  $\Phi$ we have 
\[\om^{1}_{A/B} = \langle \hd_{A/B} a_{1},\, \hd_{A/B}a_{2}, \ldots,\, \hd_{A/B}a_{r},\, \hd_{A/B}a_{r+1},\, \ldots,\, \hd_{A/B}a_{r+s}   \rangle\]
and from the exactness of  (\ref{segsucfunafdis}) it holds that 
\[
\om^{1}_{A/B} \cong \frac{\om^{1}_{B\{\mathbf{T}\}[[\mathbf{Z}]]/B}\otimes_{B\{\mathbf{T}\}[[\mathbf{Z}]]} A}{\langle \hd P_{1}\otimes1,  \ldots, \hd P_{k}\otimes1\rangle}
\]
or, equivalently, since the  functor $(-)^\tr$ is exact on $\coh(\FX)$, 
\[\om^{1}_{\FX/\FY} \cong \frac{\om^{1}_{\BD^{s}_{\BA^{r}_{\FY}}/\FY} \otimes_{\CO_{\BD^{s}_{\BA^{r}_{\FY}}}} \CO_{\FX}}{\langle \hd P_{1}\otimes1, \ldots, \hd P_{k}\otimes1\rangle^{\tr}}.\]
\end{parraf}

\section{Differentials and infinitesimal lifting properties} \label{sec3}

Next we study some characterizations  for  a smooth, unramified and \'etale morphism between locally noetherian formal schemes. Above all, we will focus on the properties related to the  module of  differentials  $\om^{1}_{\FX/\FY}$. We highlight the importance of the Jacobian Criterion for affine formal schemes (Corollary \ref{critjacob}) which allows us to determine when a closed formal subscheme of a smooth formal scheme is smooth rendering Zariski's Jacobian Criterion for topological rings (\emph{cf.} \cite[(\textbf{0}, 22.6.1)]{EGA41} into the present context.

\begin{propo} \label{condinflocal}
Let  $f:\FX \to \FY$ a morphism in $\sfn$.
\begin{enumerate}
\item
Given  $\{\FU_{\alpha}\}_{\alpha \in L}$ an open covering of  $\FX$,  $ f$  is smooth  (unramified or \'etale)  if, and only if,  for all $\alpha \in L$, $f |_{\FU_{\alpha}}: \FU_{\alpha} \to \FY$ is smooth  (unramified or \'etale, respectively).
\item
If $\{\FV_{\lambda}\}_{\lambda\in J}$ is an open  covering of  $\FY$, $ f$  is smooth  (unramified or \'etale)  if, and only if,  for all $\lambda \in J$, $f |_{f^{-1}(\FV_{\lambda})}: f^{-1}(\FV_{\lambda} ) \to \FV_{\lambda} $ is smooth  (unramified or \'etale, respectively).
\end{enumerate}
\end{propo}

\begin{proof}
This may be proved similarly as the case of usual schemes \cite[(17.1.6)]{EGA44} having in mind Propositions \ref{infsorit} and \ref{propibasinf}.
\end{proof}

\begin{cor} \label{corcondinflocal}
The results \ref{propibasinf},  \ref{corbasinf} and  \ref{condinflocal} are true if we replace the infinitesimal lifting properties by their adic counterparts.
\end{cor}
\begin{proof}
It is straightforward  from this results in view of  \cite[(10.12.1)]{EGA1} and Proposition \ref{mtf}. 
\end{proof}

\begin{rem}
It follows from the two last results that in the local study of the infinitesimal lifting properties  (with or without the adic hypothesis) over locally noetherian formal schemes, we can restrict to $\sfna$. 
\end{rem}

\begin{parraf}
We say  that $f:\FX \to \FY$ in $\sfn$ is \emph{smooth (unramified or \'etale) at $x \in \FX$ }if there exists an open subset $\FU \subset \FX$ with $x \in \FU$  such that  $f|_{\FU}$ is smooth  (unramified or \'etale, respectively).

By Proposition \ref{condinflocal} it holds that $f$ is smooth  (unramified or \'etale) if, and only if,  $f$ is smooth  (unramified or \'etale, respectively) at $x \in \FX$, $\forall x \in \FX$.
Observe that the  set of  points $x \in \FX$ such that  $f$ is smooth  (unramified, or \'etale) in $x $ is an open subset of  $\FX$.

In a forthcoming paper we will show how the infinitesimal lifting conditions in $\sfn$ at a given point depend only on the local rings.
\end{parraf}

\begin{cor} \label{complpetal}
Let $X$ be in $\sch$ and $X' \subset X$ a closed subscheme. Then the  morphism of   completion of  $X$ along  $X'$, $\kappa: X_{/X'} \to X$ is \'etale.
\end{cor}
\begin{proof}
Applying Proposition \ref{condinflocal}  we may assume  that $X= \spec(A)$ and $X'= \spec(A/I)$ are in $\scha$ with $A$ a noetherian ring and $I \subset A$ an ideal. Then, it follows from  Example \ref{exconforinf}.(\ref{compformetale}).
\end{proof}

\begin{propo} \label{complpe}
Given $f:X \to Y$ in $\sch$, let $X' \subset X$ and $Y' \subset Y$ be closed subschemes such that $f(X')\subset Y'$. 
\begin{enumerate}
\item
If $f$ is smooth  (unramified or \'etale) then $\hf: X_{/ X'} \to Y_{/ Y'}$ is smooth  (unramified or \'etale, respectively).
\item
If moreover $X' = f^{-1}(Y')$  then $\hf: X_{/ X'} \to Y_{/ Y'}$ is adic smooth  (adic unramified or adic \'etale, respectively).

\end{enumerate}
\end{propo}
\begin{proof}
Let us consider the commutative diagram of locally noetherian formal schemes
\begin{diagram}[height=2em,w=2em,labelstyle=\scriptstyle]
X             &	 \rTto^{f}  & Y\\
\uTto^{ \kappa} &	            & \uTto^{ \kappa}\\
X_{/ X'}    &	\rTto^{\hf} & Y_{/ Y'} \\
\end{diagram}
where the vertical arrows are morphisms of completion that, being slightly imprecise, we denote both by $\kappa$. Let us prove (1). If $f$ is smooth  (unramified or \'etale) by the last corollary and Proposition \ref{infsorit}.(\ref{infsorit1}) we have that $f \circ \kappa = \kappa \circ \hf$  is also smooth  (unramified or \'etale). Since $\kappa$ is \'etale from Proposition \ref{propibasinf} we deduce that $\hf$ is smooth  (unramified or \'etale, respectively). Assertion (2) is consequence of  (1)  and \cite[(10.13.6)]{EGA1}.
\end{proof}


\begin{propo} \label{modifcero}
Let $f:\FX \to \FY$ be a morphism in $\sfn$ of  pseudo finite type. The  morphism $f$ is unramified  if, and only if,  $\om^{1}_{\FX/\FY}=0$.
\end{propo}

\begin{proof}
By Proposition \ref{condinflocal} and Definition \ref{tresuno} we may suppose that $f:\FX \to \FY$ is in $\sfna$ and therefore the  result follows from \cite[(\textbf{0}, 20.7.4)]{EGA41}.
\end{proof}

\begin{cor}
Let $f:\FX \to \FY$ and $g: \FY \to \FS$ be two  pseudo finite type morphisms in $\sfn$. Then $f$ is unramified  if, and only if,  the  morphism of  $\CO_{\FX}$-Modules $f^{*}(\om^{1}_{\FY/\FS}) \to \om^{1}_{\FX/\FS}$ is surjective.
\end{cor}

\begin{proof}
Use the last proposition and the First Fundamental Exact Sequence  (\ref{sef1}) associated to  the morphisms $\FX \xto{f} \FY \xto{g} \FS$.
\end{proof}

\begin{propo}\label{flplano}(\emph{cf.} \cite[Proposition 2.6.1]{LNS})
Let $f:\FX \to \FY$ be a smooth morphism. Then $f$ is flat and $\om^{1}_{\FX/\FY}$ is  a locally free $\CO_{\FX}$-module of finite rank.
\end{propo}

\begin{proof}
Since it is a local question, we may assume that $f: \FX=\spf(A) \to \FY=\spf(B)$ is in $\sfna$ where $\phi\colon B \to A$ is a topological $B$-algebra that is formally smooth (see (\ref{condinfalg})). As for proving the flatness it suffices to show  that for all maximal ideals  $\fp \subset A$, $A_{\fp}$ is a flat $B_{\fq}$-module with  $\fq=\phi^{-1}(\fp)$. Fix $\fp \subset A$ a prime ideal. By  \cite[(\textbf{0}, 19.3.5.(iv))]{EGA41}  it holds that $A_{\fp}$ is a formally smooth $B_{\fq}$-algebra  for  the adic topologies and applying \cite[(\textbf{0}, 19.3.8)]{EGA41} there results that $A_{\fp}$ is a formally smooth $B_{\fq}$-algebra  for the topologies given by the maximal  ideals. Then, by \cite[(\textbf{0}, 19.7.1)]{EGA41} we have that $A_{\fp}$ is a flat $B_{\fq}$-module.

Let now $J \subset A$  be an ideal of  definition of $A$. The $A/J$-module $\om^{1}_{A/B}\otimes_A A/J$ is projective. Indeed, given an exact sequence $L \overset{g}\to M \to 0$ of $A/J$-modules, the sequence
\[
\Hom_{A/J}(\om^{1}_{A/B}\otimes_A A/J, L) \to 
\Hom_{A/J}(\om^{1}_{A/B}\otimes_A A/J, M) \to 0
\]
is exact by the identity of functors on $A/J$-modules
\[
\Hom_{A/J}(\om^{1}_{A/B}\otimes_A A/J, -) =
\Hom_{A}(\om^{1}_{A/B}, -)                =
\Der_B(A,-)
\]
and an argument analogous to the proof of \cite[Theorem 28.5]{ma2}.
But $\om^{1}_{A/B}$ is a  finite type $A$-module (see Proposition \ref{finitmodif}) by \cite[(\textbf{0}, 7.2.10)]{EGA1}, then it follows that $\om^{1}_{A/B}$ is a projective $A$-module. The result is now a consequence of \cite[(10.10.8.6)]{EGA1}.
\end{proof}

\begin{propo} \label{propslisoimpllocalrot}
Let $f:\FX \to \FY$ be a smooth morphism  in $\sfn$. For all pseudo finite type morphism $\FY \to \FS$  in $\sfn$ the sequence of  coherent $\CO_{\FX}$-modules 
\[
0 \to f^{*}\om^{1}_{\FY/\FS} \xto{\Phi} \om^{1}_{\FX/\FS} \xto{\Psi} \om^{1}_{\FX/\FY} \to 0
\]
defined in Proposition \ref{primersef} is exact and locally split.
\end{propo}

\begin{proof}
It is a local question, and follows from \cite[Lemma 2.5.2]{LNS} that is based on \cite[(\textbf{0}, 20.7.17.3) and (\textbf{0}, 20.7.18)]{EGA41}.
\end{proof}

\begin{cor} \label{corisomodifetal}
Let $f:\FX \to \FY$ be an \'etale morphism  in $\sfn$. For all  pseudo finite type morphism $\FY \to \FS$  in $\sfn$ it holds that
\[
f^{*}\om^{1}_{\FY/\FS} \cong \om^{1}_{\FX/\FS}
\]
\end{cor}

\begin{proof}
It is a  consequence of the last result   and of  Proposition \ref{modifcero}.
\end{proof}

\begin{parraf} \label{formallyescindida}
Given $A$ a $J$-preadic ring, let $\HA$ be the completion  of  $A$ for  the $J$-adic topology and $A_{n}=A/J^{n+1}$, for all $n \in \NN$. Take $M'',\, M'$ and $M$ $A$-modules, denote by $\widehat{M''},\,  \widehat{M' },\,\widehat{M}$ their completions for the  $J$-adic topology and let  $\widehat{M''} \xto{u} \widehat{M' } \xto{v}\widehat{M}$ be
a sequence of  $\HA$-modules. It holds that
\begin{enumerate}
\item  \label{formallyescindida1} If $0 \to \widehat{M''} \xto{u} \widehat{M' }\xto{v} \widehat{M} \to 0$ is a  split  exact sequence   of  $\HA$-modules then, for all $n \in \NN$ 
\[
0 \to M''\otimes_{A} A_{n} \xto{u_{n}} M'\otimes_{A} A_{n} \xto{v_{n}} M\otimes_{A} A_{n} \to 0
\] is a split exact sequence.

\item  \label{formallyescindida2} Reciprocally, if  $M \otimes_{A}A_{n}$ is a projective $A_{n}$-module and \[0 \to M''\otimes_{A} A_{n} \xto{u_{n}} M'\otimes_{A} A_{n} \xto{v_{n}} M\otimes_{A} A_{n} \to 0\] is a split exact sequence  of  $A_{n}$-modules, for all $n \in \NN$,  then
\begin{equation} \label{secformlinv}
0 \to \widehat{M''} \to \widehat{M' }\to \widehat{M} \to 0
\end{equation}
is a split exact sequence of $\HA$-modules. 
\end{enumerate}
Assertion (1) is immediate. In order to prove (2), for all $n \in \NN$ we have the following commutative diagrams:
\begin{diagram}[height=2em,w=1.8em,labelstyle=\scriptstyle,midshaft]
 0 &\rTto&  M'' \otimes_{A} A_{n+1}& \rTto^{u_{n+1}} &M' \otimes_{A}  A_{n+1} & \rTto^{v_{n+1}} &M \otimes_{A} A_{n+1} &  \to &0 \\
 &	&\dTonto^{f_{n}}& & \dTonto^{g_{n}}&	&\dTonto^{h_{n}}& & \\
 0 &\rTto&  M'' \otimes_{A} A_{n}& \rTto^{u_{n}} &M' \otimes_{A}  A_{n} & \rTto^{v_{n}} &M \otimes_{A} A_{n} &  \to &0 \\
\end{diagram}
where the rows are split exact sequences and the vertical maps are the canonical ones. Applying inverse limit we have that the sequence (\ref{secformlinv}) is exact. Let us show that it splits. By hypothesis, for all $n \in \NN$ there exists $t_{n}: M \otimes_{A} A_{n} \to M'\otimes_{A}  A_{n}$ such that  $v_{n} \circ t_{n}= 1$. From  $\{t_{n}\}_{n \in \NN}$ we are going to  define a family of  morphisms $\{t'_{n}: M \otimes_{A} A_{n} \to M' \otimes_{A}  A_{n}\}_{n \in \NN}$ such that 
\begin{equation} \label{conmutatividad escindidas}
v_{n} \circ t'_{n} =1 \qquad g_{n} \circ t'_{n+1} = t'_{n} \circ h_{n}
\end{equation}
for all $n \in \NN$. For  $k=0$ put $t'_{0}:=t_{0}$. Suppose that we have constructed  $t'_{k}$ verifying (\ref{conmutatividad escindidas}) for all $k \le n$ and let us define $t'_{n+1}$. If  $w_{n}:=g_{n} \circ t_{n+1} - t'_{n}\circ h_{n}$ then $v_{n} \circ w_{n} = 0$ and therefore, $\Img w_{n} \subset  \Ker v_{n}= \Img u_{n}$. Since $M \otimes_{A}A_{n+1}$ is a projective $A_{n+1}$-module, there exists $\theta_{n+1}: M\otimes_{A}  A_{n+1} \to u_{n+1}(  M''\otimes_{A} A_{n+1})$ such that  the  following diagram is commutative
\begin{diagram}[height=2em,w=1.8em,labelstyle=\scriptstyle,midshaft]
  &			& u_{n+1}( M''\otimes_{A} A_{n+1}) \\
  & \ruTto^{\theta_{n+1}}	& \dTonto\\
M \otimes_{A}  A_{n+1}& \rTto^{w_{n}} & u_{n}( M''\otimes_{A} A_{n}). \\
\end{diagram}
If we put $t'_{n+1}:=t_{n+1} - \theta_{n+1}$, it holds that $v_{n+1} \circ t'_{n+1}=1$ and $g_{n} \circ t'_{n+1} = t'_{n} \circ h_{n}$. The  morphism \[t':= \invlim {n \in \NN} t'_{n}\] satisfies that $v \circ t'=1$ and the sequence (\ref{secformlinv}) splits. 
\end{parraf}

\begin{propo}
Let $f:\FX \to \FY$ be a  pseudo finite type morphism in $\sfn$  and $g:\FY \to \FS$ a smooth morphism  in $\sfn$. The following conditions are equivalent:
\begin{enumerate}
\item
$f$ is smooth 
\item
$g \circ f $ is smooth  and the sequence 
\[
0 \to f^{*}\om^{1}_{\FY/\FS} \to \om^{1}_{\FX/\FS} \to \om^{1}_{\FX/\FY} \to 0
\]
is exact and locally split.
\end{enumerate}
\end{propo}

\begin{proof}
The  implication (1) $\Rightarrow$ (2) is consequence of  Proposition \ref{infsorit} and of  Proposition \ref{propslisoimpllocalrot}.

As for  proving that (2) $\Rightarrow$ (1) we may suppose  that $f = \spf(\phi)\colon \FX=\spf(A) \to \FY=\spf(B)$ and $g = \spf(\psi)\colon \FY=\spf(B) \to \FS=\spf(C)$ are in $\sfna$ being $B$ and $A$ formally smooth $C$-algebras. Let us show  that $A$ is a formally smooth $B$-algebra.
Let $E$ be a discrete ring  , $I \subset E$ a square zero ideal and consider the commutative diagram of  continuous homomorphisms of  topological rings
\begin{diagram}[height=2em,w=1.8em,labelstyle=\scriptstyle,midshaft]
C    & \rTto^{\psi} & B               & \rTto^{\phi} & A\\
     &              & \dTto^{\lambda} &              & \dTto^{u}\\
     &              & E               & \rTonto^{j}  & E/I.\\
\end{diagram}
Since $A$ is a formally smooth $C$-algebra, there exists a continuous homomorphism of topological $C$-algebras $v: A \to E$ such that  $v \circ \phi \circ \psi = \lambda \circ \phi$ and $j \circ v = u$. Then by   \cite[(\textbf{0}, 20.1.1)]{EGA41}) we have that $d:= \lambda - v \circ \phi \in \Dercont_{C}(B,E)$. 
From the hypothesis and  considering the equivalence of categories \ref{propitrian}.\ref{equivtrian} we have that the sequence of  finite type $A$-modules  
\[
0 \to \om^{1}_{B/C} \otimes_{B}A \to \om^{1}_{A/C}  \to \om^{1}_{A/B}   \to 0
\]
is exact and split. Besides, since the  morphism $v$ is continuous and $E$ is discrete there exists $n \in \NN$ such that  $E$ is an $A/J^{n+1}$-module.  Therefore the induced map  $\Hom_{A}(\om_{A/C}, E) \to \Hom_{B}(\om_{B/C}, E)$ is surjective and applying \cite[(\textbf{0}, 20.4.8.2)]{EGA41} we have that  the map 
\[
\Der_{C}(A,E) \to \Der_{C}(B,E)
\]
is surjective too. It follows  that there exists $d' \in \Der_{C} (A,E)$ such that  $d' \circ \phi =d$. If we put $v':= v+d'$, we have that $v' \circ \phi = \lambda$ and $j \circ v'  = u$. Therefore,  $A$ is a formally smooth $B$-algebra.
\end{proof}

\begin{cor} \label{modifisopet}
Let $f:\FX \to \FY$ and $g:\FY \to \FS$ be two  pseudo finite type morphisms in $\sfn$ such that  $g \circ f$ and $g$ are smooth. Then, $f$ is \'etale  if, and only if,  $f^{*} \om^{1}_{\FY/\FS} \cong \om^{1}_{\FX/\FS}$.
\end{cor}

\begin{proof}
Follows from the last proposition and Proposition \ref{modifcero}.
\end{proof}

\begin{propo}{(Zariski Jacobian criterion for  preadic rings)} \label{critjacobanillos}
Let $B \to A$ be a continuous morphism of  preadic rings  and suppose that $A$ is a formally smooth $B$-algebra. Given an ideal $I \subset A$. Let us consider in $A' := A/I$ the topology induced by the topology of $A$. The following conditions are  equivalent:
\begin{enumerate}
\item
$A'$ is a  formally smooth $B$-algebra.
\item
Given $J \subset A$ an ideal of  definition of  $A$, define $A'_{n}:= A/(J^{n+1}+I)$. The sequence of  $A'_{n}$-modules 
\begin{equation*}
0 \to \frac{I}{I^{2}} \otimes_{A'} A'_{n} \xto{\delta_{n}} \Omega^{1}_{A/B} \otimes_{A}   A'_{n} \xto{\Phi_{n}} \Omega^{1}_{A'/B} \otimes_{A'} A'_{n}  \to 0
\end{equation*}
is exact and  split,  for all $n \in  \NN$.
\item
The sequence of  $\widehat{A'}$-modules 
\begin{equation*}
0 \to \widehat{ \frac{I}{I^{2}}} \xto{\delta} \Omega^{1}_{A/B} \tc_{A} A' \xto{\Phi} \om^{1}_{A'/B}  \to 0
\end{equation*}
is exact and split.
\end{enumerate}
\end{propo}

\begin{proof}
The fact that (1) $\dimp$ (2) follows from  \cite[(\textbf{0}, 22.6.1), (\textbf{0}, 19.1.5) and (\textbf{0}, 19.1.7)]{EGA41} and from the Second Fundamental Exact Sequence associated to the morphisms $B \to A \to A'$. Let us show that (2) $\dimp$ (3). Since $A'$ is a formally smooth $B$-algebra,  from \cite[Theorem 28.5]{ma2} we deduce that $\Omega^{1}_{A'/B} \otimes_{A'}A'_{n}$ is a projective $A'_{n}$-module, for all $n \in \NN$  and, therefore, the result follows from \ref{formallyescindida}.(\ref{formallyescindida2}).
\end{proof}

\begin{cor}{(Zariski Jacobian criterion for formal schemes)} \label{critjacob}
Let $f: \FX=\spf(A) \to \FY=\spf(B)$ be a smooth morphism   in $\sfna$ and $\FX' \inc \FX$ a  closed immersion given by an Ideal $\,\CI=I^{\tr} \subset \CO_{\FX}$. The following conditions are equivalent:
\begin{enumerate}
\item
The composed morphism  $\FX'  \inc \FX \xto{f} \FY$ is smooth.
\item
Given $\CJ \subset \CO_{\FX}$ an Ideal of  definition, if  $\CO_{X'_{n}}:= \CO_{\FX}/(\CJ^{n+1}+\CI)$, the sequence of  coherent $\CO_{X'_{n}}$-Modules
\[
\qquad \qquad 0 \to \frac{\CI}{\CI^{2}} \otimes_{\CO_{\FX'}} \CO_{X'_{n}} \xto{\delta_{n}} \om^{1}_{\FX/\FY} \otimes_{\CO_{\FX}}   \CO_{X'_{n}} \xto{\Phi_{n}} \om^{1}_{\FX'/\FY} \otimes_{\CO_{\FX'}} \CO_{X'_{n}}  \to 0
\]
is exact  and  locally split,  for all $n \in  \NN$.
\item
The sequence of coherent $\CO_{\FX'}$-Modules
\begin{equation*}
0 \to \frac{\CI}{\CI^{2}} \xto{\delta} \om^{1}_{\FX/\FY} \otimes_{\CO_{\FX}}  \CO_{\FX'} \xto{\Phi} \om^{1}_{\FX'/\FY}  \to 0
\end{equation*}
is exact and locally split.
\end{enumerate}
\end{cor}

\begin{proof}
By Proposition \ref{finitmodif} and the equivalence of categories (\ref{propitrian}.\ref{equivtrian}) it is a consequence of the last proposition.
\end{proof}

\begin{rem}
The implication (1) $\Rightarrow$ (3) is \cite[Proposition 2.6.8]{LNS}, itself a generalization of \cite[(17.2.5)]{EGA44}.
\end{rem}

\end{document}